\documentclass[a4]{amsart}

\usepackage{amssymb}
\usepackage[all]{xy}
\usepackage{amsmath, amsfonts, amsthm , amssymb} 
\usepackage{mathrsfs}
\usepackage{color}
\usepackage{hyperref}
\usepackage{enumitem}
\usepackage{array}

\input xypic 
\usepackage{tikz}
\usepackage{tikz-cd}
\usepackage{tikz-3dplot}

\usetikzlibrary{shapes.geometric, calc}

\usetikzlibrary{decorations.markings,shapes.arrows, shapes.symbols}
    
\tikzset{paint/.style={ draw=#1!50!black, fill=#1!50 },     decorate with/.style={decorate,decoration={shape backgrounds,shape=#1,shape size=1mm}}}

\newtheorem{theorem}{Theorem}[section]
\newtheorem{proposition}[theorem]{Proposition}
\newtheorem{lemma}[theorem]{Lemma}
\newtheorem{corollary}[theorem]{Corollary}

\theoremstyle{definition}
\newtheorem{definition}[theorem]{Definition}
\newtheorem{example}[theorem]{Example}
\newtheorem{remark}[theorem]{Remark}

\newcounter{bean}

\newcommand{\R}{\mathbb{R}}
\newcommand{\C}{\mathbb{C}}
\newcommand{\PP}{\mathbb{P}}
\newcommand{\Z}{\mathbb{Z}}
 
\newcommand{\G}{\mathcal{G}}

\newcolumntype{C}[1]{>{\centering\let\newline\\\arraybackslash\hspace{0pt}}m{#1}}

\begin{document}
\title{Steenrod operations for 4-dimensional toric orbifolds}  
\author{Tseleung So} 
\address{School of Mathematics, University of Southampton, Southampton, SO17 1BJ, United Kingdom} 
\email{larry.so.tl@gmail.com}

\subjclass[2010]{Primary
55S10, 
57S12; 
Secondary
57R18, 
55P15
}
\keywords{Steenrod operations, 4-dimensional toric orbifold, stable splittings, cohomological rigidity, gauge groups, spin structures}

\begin{abstract} 
We prove necessary and sufficient conditions for the existence of non-trivial Steenrod actions on the mod-$2$ cohomology of 4-dimensional toric orbifolds. As applications, the stable homotopy type and the gauge groups of a $4$-dimensional toric orbifold are determined, a partial solution to the cohomological rigidity problem for $4$-dimensional toric orbifolds is provided, and, in the smooth case, a combinatorial criterion is established for when the toric orbifold is spin.
\end{abstract}

\maketitle

\section{Introduction} 

A \emph{toric orbifold} is an orbifold of even dimension, say of dimension $2d$, equipped with a locally standard $T^d$-action. The orbit space is a $d$-dimensional simple polytope $P$, and the $T^d$-action can be described by a \emph{characteristic function} $\lambda$ on $P$. Conversely, the toric orbifold can be recovered by a combinatorial construction based on the pair $(P,\lambda)$, which will be discussed in Subsection~\ref{subsect_combinatorial construction}. We denote by $X(P,\lambda)$ the toric orbifold associated with $(P,\lambda)$.


As there is a one-one correspondence between toric orbifolds $X(P,\lambda)$ and pairs $(P,\lambda)$ of simple polytopes and characteristic functions~\cite{soumen}, a fundamental problem is to understand the relationship between the topology of $X(P,\lambda)$, the algebraic structure of its cohomology $H^*(X(P,\lambda))$ and the combinatorial data of $(P,\lambda)$. For instance, when $X(P,\lambda)$ is smooth, that is a \emph{quasi-toric manifold}, there is a ring isomorphism
\begin{equation}\label{eqn_smooth SR/J}
H^*(X(P,\lambda);\Z)\cong SR(P)/\mathcal{J}
\end{equation}
where $SR(P)$ is the \emph{Stanley-Reisner ring} of $P$ and $\mathcal{J}$ is the ideal generated by linear relations $\lambda_{i1}x_1+\cdots+\lambda_{im}x_m$ with coefficients given by $\lambda$.
When $X(P,\lambda)$ has singularities, this isomorphism still holds after both sides are tensored with rational numbers. For more details please refer to~\cite{DJ, Dan, Jur}.
The work in~\cite{FP, BSS} extended~\eqref{eqn_smooth SR/J} by showing that if $H^*(X(P,\lambda);\Z)$ is concentrated in even degrees, then there is a ring isomorphism
\[
H^*(X(P,\lambda);\Z)\cong wSR(P,\lambda)/\mathcal{J},
\]
where $wSR(P,\lambda)$ is a subring of $SR(P)$, called the \emph{weighted Stanley-Reisner ring}, and is defined using the pair $(P,\lambda)$.
In the $4$-dimensional case, 
Fu, Song and the author derived a formula expressing the cup products of $H^*(X(P,\lambda))$ in terms of the data from $(P,\lambda)$~\cite{FSS}. Notably, in the special case where $X(P,\lambda)$ is a toric surface, our formula is ``dual'' to the intersection form from classical Intersection Theory~\cite{FSS2}.

In this paper, we focus on $4$-dimensional toric orbifolds $X(P,\lambda)$, and investigate how the pairs~$(P,\lambda)$ determine the \emph{Steenrod operations}\footnote{To be more precise, we only consider mod-$2$ Steenrod operations because, if $p$ is an odd prime, then the mod-$p$ Steenrod operations $\mathcal{P}^i$ for $i>0$ are all zero in the $4$-dimensional case.} on their cohomology.

Steenrod operations are a collection of cohomology operations
\[
\{Sq^i\colon H^j(X;\Z/2)\to H^{i+j}(X;\Z/2)\mid i=0,1,2,\cdots\}
\]
characterized by a set of axioms (see Section~\ref{sect_steenrod axiom}). They were first discovered by Steenrod who used them to study the Hopf Invariant One problem and the problem of realizing polynomial algebras as the cohomology of spaces.
They are important tools in algebraic topology as they reveal structure in topological spaces that is unseen by cup products.
For example, both $S^5\vee S^3$ and $\Sigma\C\PP^2$ have trivial cohomology ring, while $Sq^2$ acts trivially on $H^*(S^5\vee S^3;\Z/2)$ and non-trivially on $H^*(\Sigma\C\PP^2;\Z/2)$, implying that these two spaces are not homotopy equivalent.

Since $X(P,\lambda)$ is simply-connected and of dimension $4$, the Steenrod operations $Sq^i$ are always trivial for $i\neq1,2$ and 
it suffices to work with $Sq^1$ and $Sq^2$. Our main theorem provides necessary and sufficient conditions for the non-triviality of these operations.

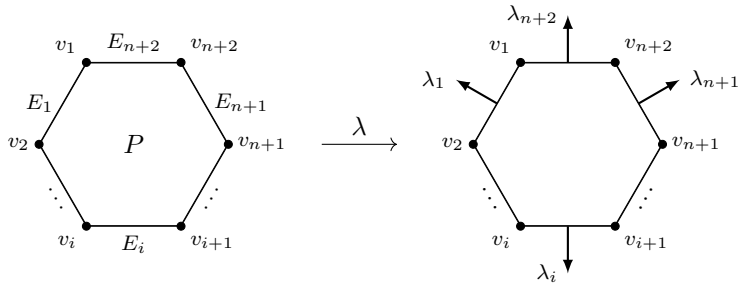
\begin{figure}
\begin{tikzpicture}
\node[opacity=0, regular polygon, regular polygon sides=6, draw, minimum size = 2.5cm](m) at (0,0) {};
\coordinate (1) at (m.corner 1); \coordinate (2) at (m.corner 2); \coordinate (3) at (m.corner 3); 
\coordinate (4) at (m.corner 4); 
\coordinate (5) at (m.corner 5); 
\coordinate (6) at (m.corner 6); 

\draw (1)--(2)--(3)--(4)--(5)--(6)--cycle;




\foreach \a in {1,...,6} { 
\draw[fill] (\a) circle (1.5pt); 
}

 \node at (0,0) {$P$};

 \node[above] at (m.side 1) {\footnotesize$E_{n+2}$};
 \node[left] at (m.side 2) {\footnotesize$E_{1}$};
 \node[right] at (m.side 6) {\footnotesize$E_{n+1}$};
 \node[below] at (m.side 4) {\footnotesize$E_{i}$};
\node[left, rotate=30] at (m.side 3) {\footnotesize$\vdots$};
 \node[right, rotate=-30] at (m.side 5) {\footnotesize$\vdots$};

 \node[above right] at (m.corner 1){\footnotesize$v_{n+2}$};
 \node[above left] at (m.corner 2){\footnotesize$v_{1}$};
 \node[left] at (m.corner 3){\footnotesize$v_{2}$};
 \node[below left] at (m.corner 4){\footnotesize$v_{i}$};
 \node[below right] at (m.corner 5){\footnotesize$v_{i+1}$};
 \node[right] at (m.corner 6){\footnotesize$v_{n+1}$};

\draw[->] (2.5,0)--(3.5,0);
\node[above] at (3,0) {$\lambda$};

\begin{scope}[xshift=150]
\node[opacity=0, regular polygon, regular polygon sides=6, draw, minimum size = 2.5cm](n) at (0.5,0) {};
\coordinate (1) at (n.corner 1); 
\coordinate (2) at (n.corner 2); 
\coordinate (3) at (n.corner 3); 
\coordinate (4) at (n.corner 4); 
\coordinate (5) at (n.corner 5); 
\coordinate (6) at (n.corner 6); 

 \draw (1)--(2)--(3)--(4)--(5)--(6)--cycle;

\foreach \a in {1,...,6} { 
\draw[fill] (\a) circle (1.5pt); 
}
\end{scope}


\node[left, rotate=30] at (n.side 3) {\small$\vdots$};
 \node[right, rotate=-30] at (n.side 5) {\small$\vdots$};

 \node[above right] at (n.corner 1){\footnotesize$v_{n+2}$};
 \node[above left] at (n.corner 2){\footnotesize$v_{1}$};
 \node[left] at (n.corner 3){\footnotesize$v_{2}$};
 \node[below left] at (n.corner 4){\footnotesize$v_{i}$};
 \node[below right] at (n.corner 5){\footnotesize$v_{i+1}$};
 \node[right] at (n.corner 6){\footnotesize$v_{n+1}$};

\draw [-latex, thick] (n.side 1) -- ($(n.side 1)!1!90:(n.corner 1)$) node [left] at ($(n.side 1)!1!90:(n.corner 1)$){\footnotesize$\lambda_{n+2}$};
\draw [-latex, thick] (n.side 2) -- ($(n.side 2)!1!90:(n.corner 2)$) node [left] at ($(n.side 2)!1!90:(n.corner 2)$){\footnotesize$\lambda_1$};
\draw [-latex, thick] (n.side 4) -- ($(n.side 4)!1!90:(n.corner 4)$) node [left] at ($(n.side 4)!1!90:(n.corner 4)$){\footnotesize$\lambda_i$};
\draw [-latex, thick] (n.side 6) -- ($(n.side 6)!1!90:(n.corner 6)$) node [right] at ($(n.side 6)!1!90:(n.corner 6)$){\footnotesize$\lambda_{n+1}$};
\end{tikzpicture}
\caption{Labels of edges and vertices in polygon $P$}\label{fig_polygon}
\end{figure}

To state our results, we need to set notation.
Let $P$ be a $2$-dimensional polygon with $n+2$ edges for $n\geq1$. Denote its edges by $E_1,\ldots,E_{n+2}$, and its vertices $v_1,\ldots,v_{n+2}$ as in Figure~\ref{fig_polygon}.

A characteristic function $\lambda$ sends the edges of $P$ to primitive vectors in $\Z^2$ that satisfy Definition~\ref{dfn_char_funct}~\eqref{dfn label_char vector independ}. For $1\leq i\leq n+2$, denote the characteristic vector $\lambda(E_i)$ by $\lambda_i$, and for $1\leq i<j\leq n+2$, denote
\[
d_{ij}=\det(\lambda_i,\lambda_j)
\]
where $(\lambda_i,\lambda_j)$ is the $2\times 2$-matrix formed by $\lambda_i$ and $\lambda_j$.

Let $X(P,\lambda)$ be the associated $4$-dimensional toric orbifold.
Thanks to~\cite{Fis, Jor, KMZ}, it is known that $H^3(X(P,\lambda);\Z)$ is a cyclic group of order
\[
g=\gcd(d_{ij}\mid 1\leq i<j\leq n+2).
\]
Therefore, for each vertex of $P$, the determinant of its adjacent characteristic vectors is a multiple of $g$.
In general, these determinants rarely equal $g$. However, due to~\cite[Lemma 5.2]{FSS0} equality can be achieved after localizing at any prime $p$. This useful fact leads to the following definition.

For any positive integer $a$, write it as $a=p^ra'$ for some index $r\geq0$ and some number $a'$ coprime to $p$. Then $r$ is the \emph{$p$-component} of $a$ and denoted by $\nu_p(a)$.

\begin{definition}\label{dfn_2-local smooth vertex}
Let $H^3(X(P,\lambda);\Z)$ have order $g$. A vertex $v_i\in P$ is called a \emph{$p$-local smooth vertex} if the determinant of its adjacent characteristic vectors has the same $p$-component as $g$. That is
\[
\nu_p(d_{i-1,i})=\nu_p(g)
\]
where $d_{0,1}$ means $d_{1,n+2}$ for $i=1$.
\end{definition}

As $2$-local smooth vertices always exist, relabelling the vertices of $P$ if necessary, we assume that $v_{n+2}$ is $2$-local smooth.
Our main theorem describes explicitly how the determinants $d_{ij}$ of characteristic vectors control the Steenrod operations on $H^*(X(P,\lambda);\Z/2)$.



\begin{theorem}\label{thm_main thm}
Let $X(P,\lambda)$ be a 4-dimensional toric orbifold associated with an $(n+2)$-gon $P$ and a characteristic function $\lambda$.
Suppose $v_{n+2}\in P$ is a $2$-local smooth vertex.
Then
\begin{enumerate}
\item\label{thm label_Sq^1 on H^2}
$Sq^1\colon H^2(X(P,\lambda);\Z/2)\to H^3(X(P,\lambda);\Z/2)$ is non-trivial if and only if
\[
g=\gcd(d_{ij}\mid 1\leq i<j\leq n+2)\equiv0\pmod{2};
\]

\item\label{thm label_Sq^1 on H^3}
$Sq^1\colon H^3(X(P,\lambda);\Z/2)\to H^4(X(P,\lambda);\Z/2)$ is always trivial;

\item\label{thm label_Sq^2}
$Sq^2\colon H^2(X(P,\lambda);\Z/2)\to H^4(X(P,\lambda);\Z/2)$ is non-trivial if and only if
\[
\prod^{n}_{i=1}\left(1-\frac{d_{i,n+1} d_{i,n+2}}{g}\right)\equiv0\pmod{2}.
\]
\end{enumerate}
\end{theorem}

\begin{remark}
Theorem~\ref{thm_main thm}~(1) and~(2) do not require $v_{n+2}$ to be a $2$-local smooth vertex, but~(3) may fail if this condition is not satisfied.
\end{remark}

This is the first calculation of the action of the Steenrod algebra on the cohomology of a toric orbifold.
The key step in proving Theorem~\ref{thm_main thm}, particularly Statement~\eqref{thm label_Sq^2}, is establishing a reduction to the special case of a triangle, which is then proved using the cup product formula derived in~\cite{FSS}.


As a corollary, Theorem~\ref{thm_main thm}~\eqref{thm label_Sq^2} implies a sufficient condition for the triviality of the $Sq^2$-action.

\begin{corollary}\label{coro_spin sufficient}
If $H^3(X(P,\lambda);\Z)$ has $2$-torsion, then $Sq^2$ acts trivially on $H^*(X(P,\lambda);\Z/2)$.
\end{corollary}

Theorem~\ref{thm_main thm}~\eqref{thm label_Sq^2} and Corollary~\ref{coro_spin sufficient} are particularly useful since the $Sq^2$-action gives a lot of information about the topology of $X(P,\lambda)$.
Four applications are given.\\

\noindent \textbf{Stable splittings.} In~\cite{ST}, the author and Theriault explored the stable splitting of a closed orientable $4$-dimensional manifold $M$ and showed that if $H_1(M;\Z)$ has no $2$-torsion, then its suspension splits into a wedge of spheres, Moore spaces, and a copy of $\Sigma\C\PP^2$ depending on its homology groups and its $Sq^2$-action. The general case is considered in~\cite{Li}, where splittings of $\Sigma^2 M$ are given that depend on more parameters and are divided into several cases.
Theorem~\ref{thm_main thm}~\eqref{thm label_Sq^2} is used to show that the stable splittings of $4$-dimensional toric orbifolds are very simple. There are only two possible homotopy types of $\Sigma X(P,\lambda)$, corresponding to the parity of $\prod^{n}_{i=1}\left(1-\frac{d_{i,n+1} d_{i,n+2}}{g}\right)$.

For integers $i\geq 2$ and $k\geq1$, the $i$-dimensional mod-$k$ \emph{Moore space} $P^i(k)$ is the mapping cone of the degree map $k\colon S^{i-1}\to S^{i-1}$.

\begin{corollary}\label{coro_stable splitting}
Let $H^3(X(P,\lambda);\Z)$ have order $g$. If $v_{n+2}$ is a $2$-local smooth vertex of $X(P,\lambda)$ and
\[
\prod^{n}_{i=1}\left(1-\frac{d_{i,n+1} d_{i,n+2}}{g}\right)\equiv0\pmod{2},
\]
%
then there is a homotopy equivalence 
\[
\Sigma X(P,\lambda)\simeq\Sigma\C\PP^2\vee\bigvee^{n-1}_{i=1}S^3\vee P^4(g).
\]
Otherwise, there is a homotopy equivalence
\[
\Sigma X(P,\lambda)\simeq S^5\vee\bigvee^{n}_{i=1}S^3\vee P^4(g).
\]
\end{corollary}



\noindent \textbf{Cohomological rigidity.}
Masuda and Suh introduced the \emph{cohomological rigidity problem}~\cite{MS08}: if two quasitoric manifolds have isomorphic cohomology rings, are they homeomorphic or diffeomorphic to each other?
This problem has attracted the attention of geometers and topologists, and affirmative answers have been proved in several cases: 2 and 4-dimensional quasitoric manifolds, certain (generalized) Bott manifolds, and 6-dimensional quasitoric manifolds associated to Pogorelov polytopes~\cite{BEMPP,CHJ,CMS-tr,OR}. Variations of this problem have also been widely studied.

Motivated by this, Fu, Song, and the author investigated the cohomological rigidity of toric orbifolds with respect to their homotopy types in~\cite{FSS0}. We proved that two $4$-dimensional toric oribfolds are homotopy equivalent if their cohomology rings are isomorphic and have no 2-torsion.
The general case, where the cohomology may contain 2-torsion, remains an open problem.

Corollaries~\ref{coro_spin sufficient} and~\ref{coro_stable splitting} are used to show that two $4$-dimensional toric orbifolds are stably homotopy equivalent if their cohomology rings are isomorphic.

\begin{corollary}\label{cor_cohmlgy rigidity}
Let $X(P,\lambda)$ and $X(P',\lambda')$ be $4$-dimensional toric orbifolds. Suppose they have isomorphic cohomology rings, and $H^3(X(P,\lambda);\Z)\cong H^3(X(P',\lambda');\Z)$ has order $g$.
\begin{itemize}
\item If $g$ is odd, then $X(P,\lambda)\simeq X(P',\lambda')$;
\item If $g$ is even, then $\Sigma X(P,\lambda)\simeq\Sigma X(P',\lambda')$.
\end{itemize}
\end{corollary}

\noindent \textbf{Gauge groups.} 
Let $X$ be a $4$-dimensional topological space and let $G$ be a simple, simply-connected compact Lie group. A principal $G$-bundle over~$X$ is classified by its second Chern class, which is an integer in this case since $H^4(BG;\Z)\cong\Z$. For any integer $k$, let $P_k$ be the principal $G$-bundle over $X$ with second Chern class $k$. The \emph{gauge group} $\G_k(X, G)$ is the topological group consisting of~$G$-equivariant automorphisms of $P_k$ that fix the base $X$.

Classifying the homotopy types of gauge groups is a longstanding problem which has been extensively studied for more than thirty years.
This has motivated topologists to develop a technique that uses stable splittings of manifolds to compute the associated gauge groups, as demonstrated in the work of~\cite{sutherland, theriault, So, huang, CS, ST}.
Here we apply the stable splitting result of Corollary~\ref{coro_stable splitting} to determine the homotopy types of gauge groups for $X(P,\lambda)$.

\begin{corollary}\label{coro_gauge group}
Let $G$ be a simple, simply-connected, compact Lie group, and let $H^3(X(P,\lambda);\Z)$ have order $g$. If $v_{n+2}$ is a $2$-local smooth vertex of $X(P,\lambda)$, $n$ is greater than $1$, and
\[
\prod^{n}_{i=1}\left(1-\frac{d_{i,n+1} d_{i,n+2}}{g}\right)\equiv0\pmod{2},
\]
%
then there is a homotopy equivalence 
\[
\G_k(X(P,\lambda),G)\simeq\G_k(\C\PP^2,G)\times\prod^{n-1}_{i=1}\Omega^2G\times\Omega^3G\{g\}
\]
where $\Omega^3G\{g\}=Map^*(P^4(g),G)$. Otherwise, there is a homotopy equivalence 
\[
\G_k(X(P,\lambda),G)\simeq\G_k(S^4,G)\times\prod^n_{i=1}\Omega^2G\times\Omega^3G\{g\}.
\]
\end{corollary}

The homotopy theory of $\G_k(S^4,G)$ and $\G_k(\C\PP^2,G)$ has been extensively studied (see~\cite{HK06,HKST18, KK18, KTT17, kono91, ST19, theriault10b,  theriault12, theriault15, theriault17}). By combining these results with Corollary~\ref{coro_gauge group}, one can classify the homotopy types of $\G_k(X(P,\lambda),G)$ for a range of Lie groups~$G$. For example, in the case where $G=SU(2)$, there is a homotopy equivalence $$\G_k(X(P,\lambda),G)\simeq \G_{k'}(X(P,\lambda),G)$$
for any integers $k$ and $k'$ if and only if
\begin{itemize}
    \item $\gcd(k,12)=\gcd(k',12)$ when $\prod^{n}_{i=1}\left(1-\frac{d_{i,n+1} d_{i,n+2}}{g}\right)$ is odd;
    \item $\gcd(k,6)=\gcd(k',6)$ when $\prod^{n}_{i=1}\left(1-\frac{d_{i,n+1} d_{i,n+2}}{g}\right)$ is even.
\end{itemize}

\noindent \textbf{Spin structures.} 
We characterize the relationship between the pair $(P,\lambda)$, the spin structure, and the second Stiefel-Whitney class when $X(P,\lambda)$ is a smooth manifold. In this case, $X(P,\lambda)$ is a \emph{quasi-toric manifold}, $H^3(X(P,\lambda);\Z)$ is trivial, and it satisfies
\[
d_{1,n+2}=\pm1,\quad\text{and}\quad
d_{i-1,i}=\pm1\quad\text{for }2\leq i\leq n+2.
\]
Therefore every vertex in $P$ is $2$-local smooth.

\begin{corollary}\label{coro_spin criteria}
Let $X(P,\lambda)$ be a 4-dimensional quasi-toric manifold. Then the following are equivalent:
\begin{enumerate}
\item
$X(P,\lambda)$ is spin;
\item
the second Stiefel-Whitney class $w_2(X(P,\lambda))$ is zero;
\item\label{coro label_spin}
we have
\[
\prod^{n}_{i=1}(1-d_{i,n+1} d_{i,n+2})\equiv1\pmod{2}.
\]
\end{enumerate}
\end{corollary}

Corollary~\ref{coro_spin criteria} recovers the results of~\cite[Corollary 6.7(ii)]{DJ} and~\cite[Proposition 0.1]{shintaro} in the $4$-dimensional case, which state that a quasi-toric manifold $X(P,\lambda)$ is spin if and only if each characteristic vector $\lambda_i=(a_i,b_i)$ satisfies
\[
a_i+b_i\equiv 1\pmod{2},
\]
that is, one of $a_i,b_i$ is even and the other one is odd.
This is equivalent to Condition~\eqref{coro label_spin} in Corollary~\ref{coro_spin criteria}.
Since $X(P,\lambda)$ is smooth, we may assume $\lambda_{n+1}=(1,0)$ and $\lambda_{n+2}=(0,1)$
after basis change in $\Z^2$. The condition $\prod^{n}_{i=1}(1-d_{i,n+1} d_{i,n+2})\equiv1\pmod{2}$
is equivalent to the condition that for all $1\leq i\leq n$
\[
a_ib_i\equiv0\pmod{2}
\]
that is, at least one of $a_i,b_i$ is even. However, $a_i$ and $b_i$ are coprime, so one of them is even and the other is odd.

The paper is organized as follows.
Section 2 introduces the notation, and reviews the combinatorial definition of toric orbifolds, key facts about their cohomology, and the $q$-CW construction from~\cite{BSS}.
Section 3 recalls the definitions of degenerate toric spaces and edge contractions introduced in~\cite{FSS}.
In Section 4 we prove preparatory lemmas, while in Section 5 we prove Theorem~\ref{thm_main thm} and the corollaries.

\section*{Acknowledgment}
The author thanks Stephen Theriault for hosting his visitor status at the University of Southampton, proofreading the draft of this paper, and giving a lot of useful feedbacks. He also thanks Xin Fu for commenting on the paper and pointing to him the work of~\cite{DJ,shintaro} on the criterion for the spin structures of quasi-toric manifolds.

\section{Preliminary facts about toric orbifolds}
\subsection{Combintorial construction of toric orbifolds}\label{subsect_combinatorial construction}
Toric orbifolds can be defined by either the combinatorial definition or the orbifold definition, which are equivalent~\cite{soumen}. For our purpose, we only present the combinatorial definition here.

A $d$-dimensional polytope~$P$ is \emph{simple} if there are exactly $d$ facets intersecting at each vertex of~$P$.
Let $\mathcal{F}(P)=\{E_1,\ldots,E_m\}$ be the set of all facets in $P$.

\begin{definition}\label{dfn_char_funct}
A \emph{characteristic function} $\lambda\colon\mathcal{F}(P)\to\Z^d$ is a map such that
\begin{enumerate}
\item
for $1\leq i\leq m$, the image $\lambda_i=\lambda(E_i)$ is a primitive vector in $\Z^d$, called the \emph{characteristic vector} of $E_i$; 
\item\label{dfn label_char vector independ}
if $E_{i_1}\cap \cdots\cap E_{i_\ell} \neq \emptyset$, then $\{\lambda_{i_1}, \dots, \lambda_{i_\ell}\}$ is linearly independent.
\end{enumerate}
\end{definition}

Let $E$ be a face of codimension-$\ell$ in $P$. As $P$ is simple, $E$ can be written as the intersection
\[
E=E_{i_1}\cap \dots \cap E_{i_\ell}
\]
of $\ell$ facets $E_{i_1},\ldots, E_{i_\ell}$, unique up to permutation. We define 
\[
T_E= \{\exp\big(t_1\lambda_{i_1}+\cdots+t_{\ell}\lambda_{i_{\ell}}\big)\mid t_1 , \ldots,t_{\ell}\in\R\}\subset T^d
\]
to be the rank-$\ell$ subtorus of $T^d$ generated by $\{\lambda_{i_1}, \dots, \lambda_{i_{\ell}}\}$.
By convention, $T_E$ is the trivial group when $E=P$.
Given a point $x\in P$, we denote by $E(x)$ the minimal face of $P$ that contains $x$ in its relative interior.

\begin{definition}\label{dfn_toric_orb}
Let $P$ be a $d$-dimensional simple polytope and let $\lambda$ be a characteristic function on $P$. Then the associated \emph{toric orbifold} $X(P,\lambda)$ is the $T^d$-space
\[
X(P, \lambda) = P \times T^d /_\sim
\]
with the equivalence relation $\sim$ and the $T^d$-action defined as follows:
\begin{itemize}
\item
$(x,t)\sim (y,s)$ if and only if $x=y$ and $t^{-1}s\in T_{E(x)}$;
\item
the $T^d$-action is defined to be
$(g,[x,t])\mapsto[x,g\cdot t]$ for any $g\in T^d$, where $[x,t]$ is the equivalence class of $(x,t)\in P\times T^d$ with respect to $\sim$.
\end{itemize}
Moreover, the \emph{orbit map} $\pi\colon X(P,\lambda)\to P$ is given by $\pi([x,t])=x$.
\end{definition}

\subsection{The cohomology of 4-dimensional toric orbifolds}\label{subsect_cohmlgy of toric orb}
The rest of the paper will focus on~4-dimensional toric orbifolds. In this case, $P$ is a $2$-dimensional polygon and facets are its edges. Unless specified otherwise, we always assume $P$ to have $n+2$ edges for $n\geq1$, and denote its edges by $E_1,\ldots,E_{n+2}$, and its vertices by $v_1,\ldots,v_{n+2}$ as in Figure~\ref{fig_polygon}.
Let $\lambda$ be a characteristic function on~$P$, and write $d_{ij}=\det(\lambda_i,\lambda_j)$ for $1\leq i<j\leq n+2$. By Definition~\ref{dfn_char_funct}~\eqref{dfn label_char vector independ}, for each vertex $v_i$ the characteristic vectors of its adjacent edges are linearly independent, or equivalently
\[
d_{1,n+2}\neq0
\quad\text{and}\quad
d_{i-1,i}\neq0\quad\text{for }2\leq i\leq n+2.
\]

The following proposition summarizes the work in~\cite{Fis, FSS0, FSS, Jor, KMZ}, which describes the group and the ring structures of the cohomology of $4$-dimensional toric orbifolds.

\begin{proposition}\label{lemma_cohomology results}
Let $X(P,\lambda)$ be a 4-dimensional toric orbifold associated with an $(n+2)$-gon $P$ and a characteristic function $\lambda$.
\begin{enumerate}[label=(\alph*)]
\item\label{lemma label_cohmlgy of toric orb}
The cohomology groups of $X(P,\lambda)$ are given as follows
\[
\begin{tabular}{C{2.5cm}|C{1cm}|C{1cm}|C{1cm}|C{1cm}|C{1cm}|C{1cm}}
$i$	&$0$	&$1$	&$2$	&$3$	&$4$	&$\geq 5$\\
\hline
$H^i(X(P,\lambda);\Z)$	&$\Z$	&$0$	&$\Z^{n}$	&$\Z/g$	&$\Z$	&$0$
\end{tabular}
\]
where $g=\gcd(d_{ij}\mid 1\leq i<j\leq n+2)$;

\item\label{lemma label_existence of 2-local smooth vert}
If $H^3(X(P,\lambda);\Z)$ has order $g$, then for any prime $p$ there exists an $i\in\{1,\ldots,n+2\}$ such that
$d_{i-1,i}$ and $g$ have the same $p$-components;

\item\label{lemma label_FSS cup prod}
If  $H^3(X(P,\lambda);\Z)$ is trivial and there exists a vertex\footnote{In~\cite{FSS} such a vertex is called a \emph{smooth vertex}. 
Note that smooth vertices do not always exist.} $v_{i}\in P$ satisfying $d_{i-1,i}=\pm1$, then there exists a set of cohomology classes
$\{u_1,\ldots,u_{n}, v\}$
such that $u_1,\ldots,u_{n}$ form a basis for~$H^2(X(P,\lambda);\Z)$, $v$ generates $H^4(X(P,\lambda);\Z)$, and
\[
u_i\cup u_j=-d_{i,n+2}d_{j,n+1}v
\]
for $1\leq i\leq j\leq n$;

\item\label{lemma label_FSS cup prod triangle}
If $P=\triangle$ is a triangle, then there exist a set of cohomology classes $\{u, v\}$ such that $u$ generates $H^2(X(\triangle,\lambda);\Z)$, $v$ generates $H^4(X(\triangle,\lambda);\Z)$, and
\[
u\cup u=-\frac{d_{12}d_{23}d_{13}}{g^2}v
\]
where $g=\gcd(d_{12},d_{23},d_{13})$.
\end{enumerate}

\end{proposition}

\begin{proof}
Statement~\ref{lemma label_cohmlgy of toric orb} is obtained from~\cite[Theorem 2.5.5]{Jor}, \cite[Theorem 2.3]{Fis}, and~\cite[Corollary 5.1]{KMZ}. Statements~\ref{lemma label_existence of 2-local smooth vert} is~\cite[Lemma 5.2]{FSS0}. Statements~\ref{lemma label_FSS cup prod} and~\ref{lemma label_FSS cup prod triangle} are Theorems 1.1 and~1.2 in~\cite{FSS} respectively.
\end{proof}

\begin{remark}
The sets $\{u_1,\ldots,u_n;v\}$ and $\{u;v\}$ in Proposition~\ref{lemma_cohomology results}~\ref{lemma label_FSS cup prod} and~\ref{lemma label_FSS cup prod triangle} are called \emph{cellular bases}.
In~\cite[Theorem 1.3]{FSS} we proved a generalization of statement~\ref{lemma label_FSS cup prod} and computed the cup products in $H^*(X(P,\lambda);R)$, where $X(P,\lambda)$ has no restrictions while $R$ is a PID satisfying certain conditions.
\end{remark}

Recall the definition of $p$-local smooth vertices in Definition~\ref{dfn_2-local smooth vertex}. Proposition~\ref{lemma_cohomology results}~\ref{lemma label_existence of 2-local smooth vert} implies that for any $X(P,\lambda)$ there always exist $p$-local smooth vertices.

\subsection{$q$-CW structures of toric orbifolds}\label{subsect_q-CW structure}

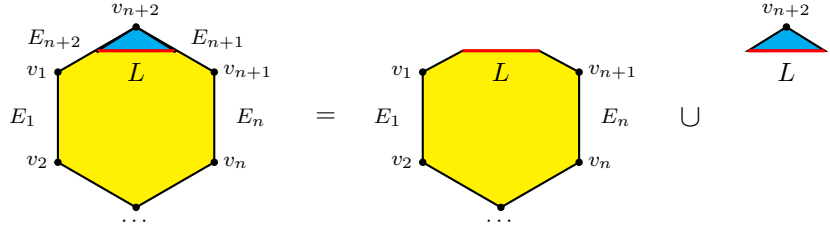
\begin{figure}
\begin{tikzpicture}[scale=0.6]

\draw[fill=yellow, thick](30:2)--(90:2)--(150:2)--(210:2)--(270:2)--(330:2)--cycle;
\draw[fill=cyan, thick] (60:1.7)--(90:2)--(120:1.7)--cycle;

\foreach \vangle in {30, 90,150,210,270,330}
\draw[fill] (\vangle:2) circle (2pt); 

\node[right] at (30:2) {\footnotesize$v_{n+1}$};
\node[above] at (90:2) {\footnotesize$v_{n+2}$};
\node[left] at (150:2) {\footnotesize$v_1$};
\node[left] at (210:2) {\footnotesize$v_2$};
\node[below] at (270:2) {\footnotesize$\cdots$};
\node[right] at (330:2) {\footnotesize$v_{n}$};

\node[right] at (60:2) {\footnotesize$E_{n+1}$};
\node[left] at (120:2) {\footnotesize$E_{n+2}$};
\node[left] at (180:2) {\footnotesize$E_1$};
\node[right] at (0:2) {\footnotesize$E_{n}$};
\draw[fill] (60:1.7) circle (1pt); 
\draw[fill] (120:1.7) circle (1pt); 
\draw[red, very thick] (60:1.7)--(120:1.7);
\node at (90:1) {$L$};
\node at (0:4.2) {$=$};

\begin{scope}[xshift=230]
\draw[fill=yellow, thick](30:2)--(60:1.7)--(120:1.7)--(150:2)--(210:2)--(270:2)--(330:2)--cycle;
\draw[red, very thick] (60:1.7)--(120:1.7);
\node at (90:1) {$L$};
\foreach \vangle in {30,150,210,270,330}
\draw[fill] (\vangle:2) circle (2pt); 

\node[right] at (30:2) {\footnotesize$v_{n+1}$};
\node[left] at (150:2) {\footnotesize$v_1$};
\node[left] at (210:2) {\footnotesize$v_2$};
\node[below] at (270:2) {\footnotesize$\cdots$};
\node[right] at (330:2) {\footnotesize$v_{n}$};

\node[left] at (180:2) {\footnotesize$E_1$};
\node[right] at (0:2) {\footnotesize$E_{n}$};

\node at (0:4.2) {\Large$\cup$};

\end{scope}

\begin{scope}[xshift=410]

\draw[fill=cyan, thick] (60:1.7)--(90:2)--(120:1.7)--cycle;

\draw[fill] (90:2) circle (2pt);
\node[above] at (90:2) {\footnotesize$v_{n+2}$};


\draw[red, very thick] (60:1.7)--(120:1.7);
\node at (90:1) {$L$};

\end{scope}
\end{tikzpicture}
\caption{$q$-CW construction of $X(P,\lambda)$}\label{fig_q CW}
\end{figure}

Following the ideas of~\cite{BSS}, we define a $q$-CW structure for a 4-dimensional toric orbifold $X(P,\lambda)$ as follows.
Choose a vertex in $P$, say $v_{n+2}$, and draw a line segment $L$ intersecting $E_{n+1}$ and $E_{n+2}$.
Then $L$ divides $P$ into two parts, a small triangle $L*\{v_{n+2}\}$, and the remaining part $P\setminus L*\{v_{n+2}\}$. 
See Figure~\ref{fig_q CW}.

Let $P'=\overline{P\setminus L*\{v_{n+2}\}}$. Then there is a decomposition
\[
X(P,\lambda)=\pi^{-1}(L*\{v_{n+2}\})\cup\pi^{-1}(P').
\]
Note that $\pi^{-1}(L*\{v_{n+2}\})$ is homeomorphic to the cone of $\pi^{-1}(L)$, so $X(P,\lambda)$ is the mapping cone of the inclusion $\pi^{-1}(L)\hookrightarrow\pi^{-1}(P')$.
One can show that $\pi^{-1}(L)$ is a lens space $S^3/G$ where $G$ is a cyclic group of order $d_{n+1,n+2}$.
Moreover, the deformation retraction $P'\to E_1\cup\cdots\cup E_{n}$ yields the homotopy equivalence
\[
\pi^{-1}(P')\simeq\pi^{-1}(E_1\cup\cdots\cup E_{n}).
\]
For $1\leq i\leq n$, let $S^2_i=\pi^{-1}(E_i)/E_i\times\{(x_0,y_0)\}$ where $(x_0,y_0)$ is the base point of $T^2$. Then $S^2_i$ is a based 2-sphere and we have
\[
\pi^{-1}(E_1\cup\cdots\cup E_{n})\cong\bigcup^{n}_{i=1}\pi^{-1}(E_i)\simeq\bigvee^{n}_{i=1}S^2_i.
\]
Let $f\colon S^3/G\to\bigvee^{n}_{i=1}S^2_i$ be the composite
\[
f\colon S^3/G\cong\pi^{-1}(L)\hookrightarrow\pi^{-1}(P')\overset{\simeq}{\longrightarrow}\bigvee^{n}_{i=1}S^2_i.
\]
Then $X(P,\lambda)$ is its mapping cone and there is a homotopy cofibration sequence
\begin{equation}\label{eqn_toric orb cofib seq}
S^3/G\overset{f}{\longrightarrow}\bigvee^{n}_{i=1}S^2_i\overset{\jmath}{\longrightarrow}X(P,\lambda)\overset{\delta}{\longrightarrow}\Sigma(S^3/G)
\end{equation}
where $\jmath$ is the composite $\bigvee^{n}_{i=1}S^2_i\simeq\pi^{-1}(\bigcup^{n}_{i=1}E_i)\hookrightarrow X(P,\lambda)$ and $\delta$ is the coboundary map.


The homotopy cofibration sequence~\eqref{eqn_toric orb cofib seq} depends on the choice of vertex $v_{n+2}$.
If $v_{n+2}$ is a $2$-local smooth vertex, then
the $2$-component of the order of $H^2(S^3/G)$ equals $\nu_2(d_{n+1,n+2})=\nu_2(g)$. This observation leads to the following lemma.

\begin{lemma}\label{lemma_decomp H^2}
If $v_{n+2}$ in the construction of~\eqref{eqn_toric orb cofib seq} is a $2$-local smooth vertex of $X(P,\lambda)$, then
\[
\jmath^*\otimes_{\Z/2}1\colon H^2(X(P,\lambda);\Z)\otimes\Z/2\overset{\cong}{\longrightarrow}H^2(\bigvee^{n}_{i=1}S^2_i;\Z)\otimes\Z/2
\]
is isomorphic and there is a group isomorphism
\begin{equation}\label{eqn_decomp H^2}
H^2(X(P,\lambda);\Z/2)\cong H^2(X(P,\lambda);\Z)\otimes{\Z/2}\oplus \delta^*(H^1(S^3/G;\Z/2)).
\end{equation}
\end{lemma}

\begin{proof}
Let $H^3(X(P,\lambda);\Z)$ have order $g$.
Apply $H^*(-;R)$ to~\eqref{eqn_toric orb cofib seq}, where $R=\Z$ or $\Z/2$, to obtain an exact sequence
\begin{equation}\label{ex seq_S^3/G-> V S^2 -> X}
0\to H^1(S^3/G)\overset{\delta^*}{\longrightarrow} H^2(X(P,\lambda))\overset{\jmath^*}{\longrightarrow}H^2(\bigvee^{n}_{i=1}S^2_i)\overset{f^*}{\longrightarrow}H^2(S^3/G)\overset{\delta^*}{\longrightarrow}H^3(X(P,\lambda))\to0.
\end{equation}
When $R=\Z$, the second $\delta^*$ is a surjection $\Z/d_{n+1,n+2}\to\Z/g$. Since $v_{n+2}$ is $2$-local smooth, we have $\nu_2(d_{n+1,n+2})=\nu_2(g)$. The exactness implies that $\jmath^*$ has cokernel $\Z/g'$ for $g'$ odd, and~\eqref{ex seq_S^3/G-> V S^2 -> X} breaks into a short exact sequence
\[
0\to H^2(X(P,\lambda))\overset{\jmath^*}{\longrightarrow}H^2(\bigvee^{n}_{i=1}S^2_i)\longrightarrow \Z/g'\to0.
\]
Since $H^2(X(P,\lambda);\Z)$ and $H^2(\bigvee^n_{i=1}S^2_i;\Z)$ are free abelian, applying $-\otimes\Z/2$ preserves the exactness and hence $\jmath^*\otimes_{\Z/2}1$ is an isomorphism.
This proves the first part of the lemma.

Take $R=\Z/2$ in~\eqref{ex seq_S^3/G-> V S^2 -> X}. By assumption, $\nu_2(d_{n+1,n+2})=\nu_2(g)$ so $H^2(S^3/G;\Z/2)\cong H^3(X(P,\lambda);\Z/2)$. The second $\delta^*$ is a surjection and hence is an isomorphism. Therefore $f^*$ is trivial and~\eqref{ex seq_S^3/G-> V S^2 -> X} breaks into a short exact sequence
\[
0\to H^1(S^3/G;\Z/2)\overset{\delta^*}{\longrightarrow} H^2(X(P,\lambda);\Z/2)\overset{\jmath^*}{\longrightarrow}H^2(\bigvee^{n}_{i=1}S^2_i;\Z/2)\to0.
\]
Since $\Z/2$ is a field, all modules are projective so we obtain
\[
H^2(X(P,\lambda);\Z/2)\cong H^2(\bigvee^{n}_{i=1}S^2_i;\Z/2)\oplus \delta^*(H^1(S^3/G;\Z/2)).
\]
It remains to identify $H^2(\bigvee^{n}_{i=1}S^2_i;\Z/2)$ with $H^2(X(P,\lambda);\Z)\otimes{\Z/2}$, which can be done via
\[
H^2(\bigvee^{n}_{i=1}S^2_i;\Z/2)\overset{\cong}{\longrightarrow} H^2(\bigvee^{n}_{i=1}S^2_i;\Z)\otimes\Z/2\overset{\jmath^*\otimes_{\Z/2}1}{\longrightarrow} H^2(X(P,\lambda);\Z)\otimes\Z/2.
\]
The first map is due to Universal Coefficient Theorem and the second map is an isomorphism due to the first part of the lemma.
Therefore isomorphism~\eqref{eqn_decomp H^2} holds.
\end{proof}

\section{Degenerate toric spaces}

The notions of degenerate toric spaces and edge contractions will be needed in later sections, so here we give a brief review.
In~\cite{FSS}, we defined degenerate toric spaces as a generalization of 4-dimensional toric orbifolds and constructed a family of maps, called \emph{toric morphisms}. Edge contractions are a special type of toric morphisms.

\begin{definition}\label{dfn_deg char fct n deg toric orb}
(see \cite[Definitions 5.1 and 5.2]{FSS})
Let $P$ be an $(n+2)$-gon with edges $E_1,\ldots,E_{n+2}$. Then a \emph{degenerate characteristic function}
\[
\lambda\colon\{E_1,\ldots,E_{n+2}\}\to\Z^2
\]
is a map that assigns to each edge $E_i$ a primitive vector $\lambda_i=\lambda(E_i)$.
Given a pair of $P$ and $\lambda$, the associated \emph{degenerate toric space} is the $T^d$-space
\[
X(P,\lambda)=P\times T^2/_{\sim_d}
\]
with the equivalence relation $\sim_d$ and the $T^d$-action defined as follows:
\begin{itemize}
\item
$(x,t)\sim_d (y,s)$ if and only if $x=y$ and
\begin{enumerate}
\item
$t^{-1}s\in T_{E(x)}$ when $x$ is not a vertex of $P$, or
\item
$t$ and $s$ are any points in $T^2$ when $x$ is a vertex of $P$;
\end{enumerate}
\item
the $T^d$-action is defined to be
$(g,[x,t]_d)\mapsto[x,g\cdot t]_d$ for any $g\in T^d$, where $[x,t]_d$ is the equivalence class of $(x,t)\in P\times T^d$ with respect to $\sim_d$.
\end{itemize}
Moreover, the \emph{orbit map} $\pi\colon X(P,\lambda)\to P$ is given by $\pi([x,t]_d)=x$.
\end{definition}

When $X(P,\lambda)$ is a degenerate toric space, we still call $\lambda_i$ the \emph{characteristic vector} of $E_i$, and write $d_{ij}=\det(\lambda_i,\lambda_j)$ for $1\leq i<j\leq n+2$.

Comparing Definitions~\ref{dfn_char_funct} and~\ref{dfn_deg char fct n deg toric orb}, we see that a degenerate characteristic function $\lambda$ becomes a characteristic function if Condition~(2) of Definition~\ref{dfn_char_funct} is satisfied, or equivalently
\[
d_{1,n+2}\neq0
\quad\text{and}\quad
d_{i,i+1}\neq0
\quad\text{for }1\leq i\leq n+1.
\]
In this case the equivalence relations $\sim$ in Definition~\ref{dfn_toric_orb} and $\sim_d$ in Definition~\ref{dfn_deg char fct n deg toric orb} coincide.

For a degenerate toric space $X(P,\lambda)$, if it satisfies $d_{n+1,n+2}\neq0$, then we can define a $q$-CW structure for $X(P,\lambda)$ as in Subsection~\ref{subsect_q-CW structure} and obtain a homotopy cofibration sequence as~\eqref{eqn_toric orb cofib seq}:
\[
S^3/G\overset{f}{\longrightarrow}\bigvee^n_{i=1}S^2_i\overset{\jmath}{\longrightarrow}X(P,\lambda)\overset{\delta}{\longrightarrow}\Sigma(S^3/G)
\]
where $G$ is a cyclic group of order $d_{n+1,n+2}$.

Here we give an example of a degenerate toric space that is not a toric orbifold.

\begin{example}\label{lemma_deg space}
Let $\triangle$ be a triangle, and let $\lambda$ be a degenerate characteristic function such that
\[
\lambda_1=\pm\lambda_3=(a,b)
\quad\text{and}\quad
\lambda_2=(c,d)
\]
where $(a,b)\neq\pm(c,d)$.
Then $d_{23}$ is not zero.
We claim that
\begin{equation}\label{eqn_deg triangle}
X(\triangle,\lambda)\simeq\Sigma(S^3/G)\vee S^2
\end{equation}
where $G$ is a cyclic group of order $d_{23}$.

In this case, the homotopy cofibration sequence~\eqref{eqn_toric orb cofib seq} becomes
\[
S^3/G\overset{f}{\longrightarrow}S^2\longrightarrow X(\triangle,\lambda).
\]
It suffices to show that $f$ is null homotopic.
Recall that $f$ is the composite
\[
f\colon S^3/G\cong\pi^{-1}(L)\hookrightarrow\pi^{-1}(\overline{\triangle\setminus L*\{v_3\}})\overset{\text{retract}}{\longrightarrow}\pi^{-1}(E_2)\simeq S^2,
\]
where $\pi^{-1}(L)=L\times T^2/\sim$ and $\pi^{-1}(E_2)=E_2\times T^2/\sim$ are quotient spaces of $[0,1]\times T^2$. Denote the quotient images of $(t, x_1, x_2)\in[0,1]\times T^2$ in $\pi^{-1}(L)$ and $\pi^{-1}(E_2)$ by $[t, x_1, x_2]_L$ and $[t, x_1, x_2]_E$ respectively. They satisfy relations
\[
\begin{array}{c}
[0, x_1, x_2]_L=[0, s^ax_1, s^bx_2]_L,\quad [1, x_1, x_2]_L=[1, (s')^c x_1, (s')^d x_2]_L\\[5pt]
[0, x_1, x_2]_E=[0, *, *]_E, \quad[t, x_1, x_2]_E=[t, (s'')^ax_1, (s'')^bx_2]_E, \quad[1, x_1, x_2]_E=[1, *, *]_E
\end{array}
\]
where $s, s', s''$ are arbitrary points in $S^1$ and $*$ is the base point of $S^1$. Using these coordinates, $f$ is given by
\[
f([t, x_1, x_2]_L)=[t, x_1,x_2]_E.
\]
Define a map $H\colon \pi^{-1}(L)\times[0,1]\to\pi^{-1}(E_2)$ by
\[
H([t, x_1, x_2]_L,\theta)=[\theta+(1-\theta)t, x_1, x_2]_E.
\]
It is well-defined since $\lambda_1$ and $\lambda_3$ are parallel.
Obviously $H(-,0)=f$ and $H(-,1)$ is a constant map, so it gives a homotopy between $f$ and a constant map and homotopy equivalence~\eqref{eqn_deg triangle} holds.\qed
\end{example}

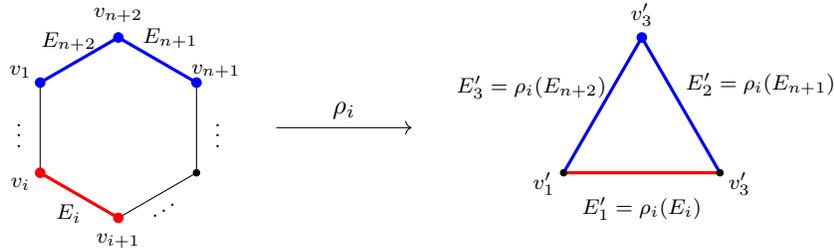
\begin{figure}
\begin{tikzpicture}[scale=0.6]

\draw (30:2)--(90:2)--(150:2)--(210:2)--(270:2)--(330:2)--cycle;

\draw[blue, very thick] (30:2)--(90:2)--(150:2);
\draw[red, very thick] (270:2)--(210:2);


\node at (30:2.5) {\footnotesize$v_{n+1}$};
\node at (90:2.5) {\footnotesize$v_{n+2}$};
\node at (150:2.5) {\footnotesize$v_1$};
\node at (210:2.5) {\footnotesize$v_i$};
\node at (270:2.5) {\footnotesize$v_{i+1}$};
\node at (330:2.5) {\text{ }};

\node at (60:2.3) {\footnotesize$E_{n+1}$};
\node at (0:2.2) {\footnotesize$\vdots$};
\node[rotate=30] at (300:2.1) {\footnotesize$\cdots$};
\node at (240:2.2) {\footnotesize$E_i$};
\node at (180:2.2) {\footnotesize$\vdots$};
\node at (120:2.2) {\footnotesize$E_{n+2}$};

\foreach \a in {30,90,...,330} { 
\draw[fill] (\a:2) circle (2pt); 
}

\draw[fill, blue] (90:2) circle (3pt); 
\draw[fill, blue] (30:2) circle (3pt); 
\draw[fill, blue] (150:2) circle (3pt); 

\draw[fill, red] (270:2) circle (3pt); 
\draw[fill, red] (210:2) circle (3pt);


\draw[->] (3.5,0)--(6.5,0);
\node[above] at (5,0) {$\rho_i$};

\begin{scope}[xshift=330, yshift=0]

\draw (90:2)--(210:2)--(330:2)--cycle;
\draw[blue, very thick] (90:2)--(210:2);
\draw[red, very thick] (210:2)--(330:2);
\draw[blue, very thick] (330:2)--(90:2);

\foreach \a in {90,210,330} { 
\draw[fill] (\a:2) circle (2pt); 
}
\node at (90:2.5) {\footnotesize$v_3'$};
\node at (210:2.5) {\footnotesize$v_1'$};
\node at (330:2.5) {\footnotesize$v_3'$};

\node at (160:2.6) {\footnotesize$E_3'=\rho_i(E_{n+2})$};
\node at (270:1.8) {\footnotesize$E_1'=\rho_i(E_i)$};
\node at (20:2.8) {\footnotesize$E_2'=\rho_i(E_{n+1})$};

\draw[fill,blue] (90:2) circle (3pt); 

\end{scope}
\end{tikzpicture}
\caption{Edge contraction $\rho_i$}\label{fig_edge contraction}
\end{figure}

Now we define edge contractions between degenerate toric spaces.
Let $P$ be a polygon with edges~$E_1,\ldots,E_{n+2}$, and let $\triangle$ be a triangle with edges $E'_1,E'_2,E'_3$.
For $1\leq i\leq n$, let $\rho_i\colon\partial P\to\partial\triangle$ be a map contracting all edges in $P$ except $E_i,E_{n+1}$ and $E_{n+2}$ (see Figure~\ref{fig_edge contraction}). Apply $Cone(-)$ to extend $\rho_i$ throughout the whole $P$ to obtain a map
\[
\rho_i\colon P\cong Cone(\partial P)\longrightarrow Cone(\partial\triangle)\cong\triangle.
\]

\begin{definition}\label{dfn_edge contraction}
(see~\cite[Definition 5.4]{FSS})
Let $P$ be an $(n+2)$-gon, and let $\lambda$ be a degenerate characteristic function on it.
For~$1\leq i\leq n$, let $\rho_i\colon P\to\triangle$ be the map given above. Then we define a degenerate characteristic function $\rho_i\lambda$ on $\triangle$ by
\[
\rho_i\lambda(E'_1)=\lambda_i,\quad
\rho_i\lambda(E'_2)=\lambda_{n+1},\quad
\rho_i\lambda(E'_3)=\lambda_{n+2},
\]
and a map $\rho_i\colon X(P,\lambda)\to X(\triangle,\rho_i\lambda)$, called an \emph{edge contraction}, by
\[
\rho_i([x, t_1,t_2]_d)=[\rho_i(x), t_1,t_2]'_d
\]
where $[x, t_1,t_2]_d$ is the equivalence class of $(x, t_1,t_2)\in P\times T^2$ in $X(P,\lambda)$ and $[\rho_i(x), t_1,t_2]'_d$ is the equivalence class of $(\rho(x), t_1,t_2)\in \triangle\times T^2$ in $X(\triangle,\rho_i\lambda)$.
\end{definition}

When $X(P,\lambda)$ is a toric orbifold, $X(\triangle,\rho_i\lambda)$ is generally not a toric orbifold but exhibits similar properties.
By Definition~\ref{dfn_edge contraction}, $\det((\rho_i\lambda)_{2},(\rho_i\lambda)_3)=\det(\lambda_{n+1},\lambda_{n+2})$.
The latter determinant is not zero, so homotopy cofibration sequence~\eqref{eqn_toric orb cofib seq} holds for $X(\triangle,\rho_i\lambda)$. Moreover, the argument in the proof of \cite[Lemma 5.13]{FSS} shows that $\rho_i$ is compatible with these homotopy cofibration sequences, yielding a homotopy commutative diagram
\begin{equation}\label{dgrm_edge contract}
\xymatrix{
S^3/G\ar[r]^-{f}\ar[d]^-{=}	&\bigvee^{n}_{j=1}S^2_j\ar[r]^-{\jmath}\ar[d]^-{\varpi_i}	&X(P,\lambda)\ar[r]^-{\delta}\ar[d]^-{\rho_i}	&\Sigma(S^3/G)\ar[d]^-{=}\\
S^3/G\ar[r]^-{f_i}	&S^2_i\ar[r]^-{\jmath_i}								&X(\triangle,\rho_i\lambda)\ar[r]^-{\delta_i}	&\Sigma(S^3/G)
}
\end{equation}
where $\varpi_i$ is the canonical pinch map.
If we further assume that $v_{n+2}\in P$ is a $2$-local smooth vertex, then $H^*(X(\triangle,\rho_i\lambda))$ reflects significant information about $H^*(X(P,\lambda))$ through $\rho^*_i$.

\begin{lemma}\label{lemma_edge contraction}
Let $X(P,\lambda)$ be a 4-dimensional toric orbifold such that $H^3(X(P,\lambda);\Z)$ has order $g$. Suppose $v_{n+2}$ is a $2$-local smooth vertex.
For $1\leq i\leq n$, let $X(\triangle,\rho_i\lambda)$ and
$\rho_i\colon X(P,\lambda)\to X(\triangle,\rho_i\lambda)$
be the degenerate toric space and the edge contraction given in Definition~\ref{dfn_edge contraction}.
\begin{enumerate}[label=(\alph*)]
\item\label{lemma label_cohmlgy of deg toric sp}
The cohomology groups of $X(\triangle,\rho_i\lambda)$ are given as follows:
\begin{equation}\label{table_deg toric space cohmlgy}
\begin{tabular}{C{2.8cm}|C{1cm}|C{1cm}|C{1cm}|C{1cm}|C{1cm}|C{1cm}}
$j$	&$0$	&$1$	&$2$	&$3$	&$4$	&$\geq 5$\\
\hline
$H^j(X(\triangle,\rho_i\lambda);\Z)$	&$\Z$	&$0$	&$\Z$	&$\Z/g_i$	&$\Z$	&$0$
\end{tabular}
\end{equation}
where $g_i$ satisfies $\nu_2(g_i)=\nu_2(g)$;

\item\label{lemma label_edge contract deg 4 isom}
The induced homomorphism
$\rho^*_i\colon H^4(X(\triangle,\rho_i\lambda);\Z/2)\to H^4(X(P,\lambda);\Z/2)$ is isomorphic;

\item\label{lemma label_edge contract deg 2 basis}
For each $i\in\{1,\ldots,n\}$ choose a generator $u_i$ of $H^2(X(\triangle,\rho_i\lambda);\Z)$.
Then 
\[
\rho^*_1(u_1)\otimes_{\Z/2}1,\ldots,\rho^*_{n}(u_{n})\otimes_{\Z/2}1
\]
form a basis for $H^2(X(P,\lambda);\Z)\otimes\Z/2$.
\end{enumerate}
\end{lemma}

\begin{proof}
First we prove~\ref{lemma label_cohmlgy of deg toric sp}. Suppose $X(\triangle,\rho_i\lambda)$ is a toric orbifold.  By Proposition~\ref{lemma_cohomology results}~\ref{lemma label_cohmlgy of toric orb}, its cohomology groups take the form shown in Table~\eqref{table_deg toric space cohmlgy}, where $g_i=\gcd(d_{i,n+1}, d_{i,n+2}, d_{n+1,n+2})$. 
It remains to show that $\nu_2(g_i)=\nu_2(g)$.
Since~$v_{n+2}$ is $2$-local smooth, we have
\[
\nu_2(d_{n+1,n+2})=\nu_2(g)=\min\{\nu_2(d_{i,j})\mid 1\leq i<j\leq n+2\}.
\]
It follows that $\nu_2(g_i)=\min\{\nu_2(d_{i,n+1}),\nu_2(d_{i,n+2}),\nu_2(d_{n+1,n+2})\}=\nu_2(d_{n+1,n+2})=\nu_2(g).$

Second, we prove~\ref{lemma label_cohmlgy of deg toric sp} when $X(\triangle,\rho_i\lambda)$ is a not toric orbifold. From the computation in Example~\ref{lemma_deg space} we know that $X(\triangle,\rho_i\lambda)\simeq\Sigma(S^3/G)\vee S^2$, where $G$ is a cyclic group of order $d_{n+1,n+2}$. Therefore its cohomology groups take the form shown in Table~\eqref{table_deg toric space cohmlgy}, where~$g_i=d_{n+1,n+2}$. The assumption implies that $\nu_2(g_i)=\nu_2(d_{n+1,n+2})=\nu_2(g)$.

Third, we prove~\ref{lemma label_edge contract deg 4 isom}. Apply $H^*(-;\Z/2)$ to~\eqref{dgrm_edge contract} to obtain a commutative diagram
\[
\xymatrix
{
0\ar[r]	&H^3(S^3/G)\ar[r]^-{\delta^*_i}\ar[d]^-{=}	&H^4(X(\triangle,\rho_i\lambda))\ar[r]\ar[d]^-{\rho^*_i}	&0\\
0\ar[r]	&H^3(S^3/G)\ar[r]^-{\delta^*}				&H^4(X(P,\lambda))\ar[r]			&0
}
\]
Since both $\delta^*_i$ and $\delta^*$ are isomorphic, so is $\rho^*_i\colon H^4(X(\triangle,\rho_i\lambda);\Z/2)\to H^4(X(P,\lambda);\Z/2)$.

Fourth, we prove~\ref{lemma label_edge contract deg 2 basis}.
Apply $H^*(-;\Z)$ to~\eqref{dgrm_edge contract} to obtain a commutative diagram
\[
\xymatrix{
0\ar[r]	&H^2(X(\triangle,\rho_i\lambda))\ar[r]^-{\jmath_i^*}\ar[d]^-{\rho^*_i}	&H^2(S^2_i)\ar[r]^-{f^*_i}\ar[d]^-{\varpi^*_i}	&H^2(S^3/G)\ar[r]^-{\delta^*_i}\ar[d]^-{=}	&H^3(X(\triangle,\rho_i\lambda))\ar[r]\ar[d]^-{\rho^*_i}	&0\\
0\ar[r]	&H^2(X(P,\lambda))\ar[r]^-{\jmath^*}	&\bigoplus^{n}_{j=1}H^2(S^2_j)\ar[r]^-{f^*}	&H^k(S^3/G)\ar[r]^-{\delta^*}	&H^3(X(P,\lambda))\ar[r]	&0
}
\]
From~\ref{lemma label_cohmlgy of deg toric sp}, we have $\nu_2(g_i)=\nu_2(d_{n+1,n+2})$. 
Consequently, the argument of Lemma~\ref{lemma_decomp H^2} also applies to~$X(\triangle,\rho_i\lambda)$ so~$\jmath^*_i\otimes_{\Z/2}1$ is an isomorphism and $\jmath^*_i(u_i)\otimes_{\Z/2}1$ is a generator of $H^2(S^2_i;\Z)\otimes\Z/2$. Since $\varpi^*_i$ is the canonical inclusion for $1\leq i\leq n$, the cohomology classes
\[
\varpi^*_1\circ\jmath^*_1(u_1)\otimes_{\Z/2}1,\ldots,\varpi^*_n\circ\jmath^*_n(u_{n})\otimes_{\Z/2}1
\]
form a basis for $\bigoplus^{n}_{j=1}H^2(S^2_j;\Z)\otimes\Z/2$. The left square implies that for $1\leq i\leq n$
\[
\varpi^*_i\circ\jmath^*_i(u_i)\otimes_{\Z/2}1=\jmath^*\circ\rho^*_i(u_i)\otimes_{\Z/2}1
\]
and $\jmath^*\otimes_{\Z/2}1$ is an isomorphism by Lemma~\ref{lemma_decomp H^2}, so
\[
\rho^*_1(u_1)\otimes_{\Z/2}1,\ldots,\rho^*_{n}(u_{n})\otimes_{\Z/2}1
\]
form a basis for $H^2(X(P,\lambda);\Z)\otimes\Z/2$.
\end{proof}

\section{Steenrod operations on the cohomology of toric orbifolds}\label{sect_steenrod axiom}

Steenrod operations are a collection of cohomological operations
\[
\{Sq^i\colon H^j(X;\Z/2)\longrightarrow H^{i+j}(X;\Z/2)\mid i=0,1,2,\cdots\}
\]
characterized by the following axioms (for reference, please see~\cite{Hat} and~\cite{MT}):
\begin{enumerate}
\item
$Sq^0$ is the identity map, $Sq^1$ is the Bockstein homomorphism $\beta$ for the exact sequence
\[
0\longrightarrow\Z/2\longrightarrow\Z/4\longrightarrow\Z/2\longrightarrow0;
\]
\item\label{axiom_Sq u = u^2}
if $x\in H^i(X;\Z/2)$ then $Sq^i(x)=x\cup x$ and $Sq^j(x)=0$ for $j>i$;
\item
if $h\colon X\to Y$ is a map, then $Sq^i(h^*(x))=h^*(Sq^i(x))$ for $x\in H^*(Y;\Z/2)$;
\item
if $\sigma\colon H^i(X)\to H^{i+1}(\Sigma X)$ is the suspension isomorphism, then $Sq^i(\sigma(x))=\sigma( Sq^i(x))$;
\item
they satisfy the Cartan Formula $Sq^i(x\cup y)=\sum_{j+k=i}Sq^j(x)\cup Sq^k(y)$;
\item
they satisfy the Adem relations
\[
Sq^iSq^j\equiv\sum^{\lfloor{i/2}\rfloor}_{k=0}\binom{j-k-1}{i-2k}Sq^{i+j-k}Sq^k\pmod{2}
\]
for all $i,j>0$ such that $i<2j$.
\end{enumerate}

In this section, we compute Steenrod operations on $H^*(X(P,\lambda);\Z/2)$, where $X(P,\lambda)$ is a 4-dimensional toric orbifold unless specified otherwise.
To begin with, let us consider homotopy cofibration sequence~\eqref{eqn_toric orb cofib seq}
\[
S^3/G\overset{f}{\longrightarrow}\bigvee^{n}_{i=1}S^2_i\longrightarrow X(P,\lambda).
\]
The lens space $S^3/G$ is a 3-dimensional orientable manifold, so it is parallelizable. Thanks to~\cite{atiyah} its double suspension has a splitting
\[
\Sigma^2(S^3/G)\simeq S^5\vee P^5(d_{n+1,n+2}),
\]
where $P^5(d_{n+1,n+2})$ is the mapping cone of the degree map 
$S^4\overset{d_{n+1,n+2}}{\longrightarrow}S^4$.
Let $f'$ be the composite
\[
f'\colon S^5\hookrightarrow S^5\vee P^5(d_{n+1,n+2})\simeq\Sigma^2(S^3/G)\overset{\Sigma^2f}{\longrightarrow}\bigvee^{n}_{i=1}S^4_i,
\]
and let $X'$ be the mapping cone of $f'$. Then there is a homotopy commutative diagram
\begin{equation}\label{dgrm_sigma^2 RP split}
\xymatrix{
S^5\ar[r]^-{f'}\ar[d]^-{\text{incl}}		&\bigvee^{n}_{i=1}S^4_i\ar[r]\ar[d]^-{=}	&X'\ar[r]\ar[d]^-{\Psi}	&S^6\ar[d]^-{\text{incl}}\\
\Sigma^2(S^3/G)\ar[r]^-{\Sigma^2 f}\ar[d]^-{\text{pinch}}	&\bigvee^{n}_{i=1}S^4_i\ar[r]\ar[d]		&\Sigma^2 X(P,\lambda)\ar[r]^-{\Sigma^2\delta}\ar[d]	&\Sigma^3(S^3/G)\ar[d]^-{\text{pinch}}\\
P^5(d_{n+1,n+2})\ar[r]	&\ast\ar[r]			&P^6(d_{n+1,n+2})\ar[r]^-{=}				&P^6(d_{n+1,n+2})
}
\end{equation}
where $\Psi$ is an induced map. The rows and columns are homotopy cofibration sequences, while the middle row is the double suspension of~\eqref{eqn_toric orb cofib seq}.
Using the third column, one can show that
\[
\Psi^*\colon H^6(\Sigma^2X(P,\lambda);\Z/2)\longrightarrow H^6(X';\Z/2)
\]
is an isomorphism.

\begin{lemma}\label{lemma_Sq^1 on H^3}
The Steenrod operation $Sq^1\colon H^3(X(P,\lambda);\Z/2)\to H^4(X(P,\lambda);\Z/2)$ is trivial.
\end{lemma}

\begin{proof}
Prove by contradiction. Assume there exists $v\in H^3(X(P,\lambda);\Z/2)$ such that $Sq^1(v)\neq0$.
Since $Sq^1$ is stable, we have $Sq^1(\sigma^2(v))\neq0$ where $\sigma$ is the suspension isomorphism.
Then
\[
\Psi^*(Sq^1(\sigma^2(v)))=Sq^1(\Psi^*(\sigma^2(v)))
\]
is not zero since $\Psi^*$ is isomorphic in degree $6$.
However, $X'$ is a CW-complex with $4$-cells and one~$6$-cell only, so $H^5(X';\Z/2)=0$ and $\Psi^*(\sigma^2(v))=0$. This leads to contradiction.
\end{proof}

Next, we turn our attention to $Sq^2\colon H^2(X(P,\lambda);\Z/2)\to H^4(X(P,\lambda);\Z/2)$. The idea is using edge contractions
\[
\rho_i\colon X(P,\lambda)\longrightarrow X(\triangle,\rho_i\lambda)
\]
defined in Definition~\ref{dfn_edge contraction} to reduce the problem to the special case of $X(\triangle,\rho_i\lambda)$. This will be done in Lemma~\ref{lemma_Sq on P = Sq on triangle}. Before proceeding, we first prove a preparatory lemma. By Lemma~\ref{lemma_decomp H^2}~\eqref{eqn_decomp H^2}
\[
H^2(X(P,\lambda);\Z/2)\cong H^2(X(P,\lambda);\Z)\otimes\Z/2\oplus\delta^*(H^1(S^3/G;\Z/2)).
\]
We show that $Sq^2$ acts trivially on $\delta^*(H^1(S^3/G;\Z/2))$.

\begin{lemma}\label{lemma_non trivial Sq is integral}
If $v_{n+2}$ is a $2$-local smooth vertex, then the restriction of $Sq^2$ to $\delta^*(H^1(S^3/G;\Z/2))$ is trivial.
Moreover, if $Sq^2$ acts on $H^*(X(P,\lambda);\Z/2)$ non-trivially, then there exists an integral cohomology class $u\in H^2(X(P,\lambda);\Z)$ such that $Sq^2(u\otimes_{\Z/2}1)\neq0$.
\end{lemma}

\begin{proof}
Let $w\in H^1(S^3/G;\Z/2)$ be the generator. To prove the first part, it suffices to show that
\[
Sq^2(\delta^*(w))=0.
\]
If not, then $Sq^2(\sigma^2(\delta^*w))\neq0$ where $\sigma$ is the suspension isomorphism.
Apply $\Psi^*$ to obtain
\[
\Psi^*(Sq^2(\sigma^2(\delta^*w))=Sq^2(\Psi^*\circ\sigma^2(\delta^*w)).
\]
Since $\Psi^*$ is isomorphic in degree $6$, $\Psi^*(Sq^2(\sigma^2(\delta^*w))\neq0$.

On the other hand, in the top right square of~\eqref{dgrm_sigma^2 RP split} $\Sigma^2\delta\circ\Psi$ factors through $S^6$, implying that
\[
\Psi^*\circ\delta^*(\sigma^2w)=(\Sigma^2\delta\circ\Psi)^*(\sigma^2w)=0
\]
as $\sigma^2w$ is of degree $3$. Therefore $Sq^2(\Psi^*\circ\delta^*(\sigma^2w))=0$, leading to contradiction. It must be the case that $Sq^2(\delta^*(w))=0$ and the restriction of $Sq^2$ to $\delta^*(H^1(S^3/G);\Z/2)$ is trivial.


The second part of the lemma follows immediately from the first part and Lemma~\ref{lemma_decomp H^2}~\eqref{eqn_decomp H^2}.
\end{proof}

\begin{lemma}\label{lemma_Sq on P = Sq on triangle}
Let $X(P,\lambda)$ be a 4-dimensional toric orbifold. Suppose $v_{n+2}$ is a $2$-local smooth vertex.
For $1\leq i\leq n$, let $\rho_i\colon P\to\triangle$ be the map contracting all edges of $P$ except $E_i,E_{n+1}$ and~$E_{n+2}$ and let~$X(\triangle,\rho_i\lambda)$ be the degenerate toric space given in Definition~\ref{dfn_edge contraction}.
Then the following are equivalent:
\begin{enumerate}[label=(\alph*)]
\item
$H^*(X(P,\lambda);\Z/2)$ has non-trivial $Sq^2$-action;

\item
at least one of $H^*(X(\triangle,\rho_i\lambda);\Z/2)$'s has non-trivial $Sq^2$-action.
\end{enumerate}
\end{lemma}

\begin{proof}
First we show $(a)\Rightarrow(b)$. If $Sq^2$ acts non-trivially on $H^*(X(P,\lambda);\Z/2)$, then Lemma~\ref{lemma_non trivial Sq is integral} implies that there is an integral cohomology class $u\in H^2(X(P,\lambda);\Z)$ such that $Sq^2(u\otimes_{\Z/2}1)\neq0$.

For $1\leq i\leq n$, let $u_i$ be a generator of $H^2(X(\triangle,\rho_i\lambda);\Z)$. By Lemma~\ref{lemma_edge contraction}~\ref{lemma label_edge contract deg 2 basis}
\[
\rho^*_1(u_1)\otimes_{\Z/2}1,\ldots,\rho^*_{n}(u_{n})\otimes_{\Z/2}1
\]
form a basis for $H^2(X(P,\lambda);\Z)\otimes\Z/2$. For convenience, write $\bar{u}_i=\rho^*_i(u_i)\otimes_{\Z/2}1$.
Then
\[
Sq^2(u\otimes_{\Z/2}1)=Sq^2\left(\sum^{n}_{i=1}a_i\bar{u}_i)\right)
=\sum^{n}_{i=1} a_i Sq^2(\bar{u}_i)
\]
for some coefficients $a_1,\ldots,a_{n}\in\Z/2$. Since $Sq^2(u\otimes_{\Z/2}1)\neq0$, at least one of $a_iSq^2(\bar{u}_i)$'s is not zero.
Therefore there exists an $i\in\{1,\ldots,n\}$ such that
\[
Sq^2(\bar{u}_i)=Sq^2(\rho^*_i(u_i)\otimes_{\Z/2}1)=\rho^*_i(Sq^2(u_i\otimes_{\Z/2}1))\neq0.
\]
By Lemma~\ref{lemma_edge contraction}~\ref{lemma label_edge contract deg 4 isom} $\rho^*_i$ is an isomorphism in degree $4$, so $Sq^2(u_i\otimes_{\Z/2}1)$ is not zero and $Sq^2$ acts non-trivially on~$H^2(X(\triangle,\rho_i\lambda);\Z/2)$.

Next we prove $(b)\Rightarrow(a)$.
Without loss of generality, assume that $H^*(X(\triangle,\rho_1\lambda);\Z/2)$ has non-trivial $Sq^2$-action. Then there exists $v\in H^*(X(\triangle,\rho_1\lambda);\Z/2)$ such that $Sq^2(v)\neq0$.
It follows that 
$Sq^2(\rho^*_1(v))=\rho^*_1(Sq^2(v))$ is not zero since $\rho^*_1$ is isomorphic in degree $4$.
Therefore $Sq^2$ acts non-trivially on $H^*(X(P,\lambda);\Z/2)$.
\end{proof}

Due to Lemma~\ref{lemma_Sq on P = Sq on triangle}, it suffices to consider the $Sq^2$ operations on the cohomology of degenerate toric spaces $X(\triangle,\rho_i\lambda)$.

\begin{lemma}\label{lemma_criteria of Sq on triangle}
Let $X(P,\lambda)$ be a $4$-dimensional toric orbifold such that $v_{n+2}$ is a $2$-local smooth vertex.
For $1\leq i\leq n$, let $X(\triangle,\rho_i\lambda)$ be the degenerate toric space given in Definition~\ref{dfn_edge contraction}.
Then the following are equivalent:
\begin{enumerate}[label=(\alph*)]
\item
$H^*(X(\triangle,\rho_i\lambda);\Z/2)$ has non-trivial $Sq^2$-action
\item
we have $\displaystyle
\frac{d_{i,n+1}\cdot d_{i,n+2}\cdot d_{n+1,n+2}}{g^2_i}\equiv1\pmod{2}$, where $g_i=\gcd(d_{i,n+1}, d_{i,n+2}, d_{n+1,n+2})$.
\end{enumerate}
\end{lemma}

\begin{proof}
For convenience, write $\tilde{\lambda}_j=\rho_i\lambda(E_j)$ for $1\leq j\leq n+2$.
By Definition~\ref{dfn_edge contraction} we have
\[
\tilde{\lambda}_1=\lambda_i,\quad
\tilde{\lambda}_2=\lambda_{n+1},\quad
\text{and}\quad
\tilde{\lambda}_3=\lambda_{n+2}.
\]
For $1\leq j<k\leq n+2$, write
$\tilde{d}_{jk}=\det(\tilde{\lambda}_j,\tilde{\lambda}_k)$. Then
\[
\tilde{d}_{12}=d_{i,n+1},\quad
\tilde{d}_{13}=d_{i,n+2},
\quad\text{and}\quad
\tilde{d}_{23}=d_{n+1,n+2}.
\]

First we prove the lemma when $X(\triangle,\rho_i\lambda)$ is a toric orbifold.
Proposition~\ref{lemma_cohomology results}~\ref{lemma label_cohmlgy of toric orb} and~\ref{lemma label_FSS cup prod triangle}  imply that the order of~$H^3(X(\triangle,\rho_i\lambda);\Z)$ equals
\[
g_i=\gcd(\tilde{d}_{12},\tilde{d}_{13},\tilde{d}_{23})=\gcd(d_{i,n+1},d_{i,n+2},d_{n+1,n+2}),
\]
and there exists a generator $u\in H^2(X(\triangle,\rho_i\lambda);\Z)$ such that
\[
u\cup u=\frac{\tilde{d}_{12}\tilde{d}_{13}\tilde{d}_{23}}{g^2_i}=\frac{d_{i,n+1}d_{i,n+2}d_{n+1,n+2}}{g^2_i}.
\]
Writing $\bar{u}=u\otimes_{\Z/2}1$ and using Axiom~\ref{axiom_Sq u = u^2} of Steenrod operations, we have
\[
Sq^2(\bar{u})=\bar{u}\cup\bar{u}=\frac{d_{i,n+1}d_{i,n+2}d_{n+1,n+2}}{g^2_i}.
\]
Therefore $Sq^2$ acts non-trivially on $H^*(X(\triangle,\rho_i\lambda);\Z/2)$ if and only if
\[
\frac{d_{i,n+1}d_{i,n+2}d_{n+1,n+2}}{g^2_i}\equiv1\pmod{2}.
\]

Now, assume $X(\triangle,\rho_i\lambda)$ is not a toric orbifold. 
In this case, at least one of $\tilde{d}_{12},\tilde{d}_{13},\tilde{d}_{23}$ is zero, so
\[
\frac{d_{i,n+1}d_{i,n+2}d_{n+1,n+2}}{g^2_i}=0.
\]
It remains to show that $H^*(X(\triangle,\rho_i\lambda);\Z/2)$ has trivial $Sq^2$-action.
Notice that $\tilde{d}_{23}=d_{n+1,n+2}$ is not zero, so either $\tilde{d}_{12}$ or $\tilde{d}_{13}$ is zero. Therefore, $\tilde{\lambda}_1$ is parallel to either $\tilde{\lambda}_2$ and $\tilde{\lambda}_3$, but the latter two vectors are not parallel. This is the case considered in Example~\ref{lemma_deg space}, where we have shown
\[
X(\triangle,\rho_i\lambda)\simeq\Sigma(S^3/G)\vee S^2
\]
with $G$ a cyclic group of order $\tilde{d}_{23}=d_{n+1,n+2}$.
As $H^*(S^2;\Z/2)$ and $H^*(\Sigma(S^3/G);\Z/2)$ have trivial~$Sq^2$-actions, so does $H^*(X(\triangle,\rho_i\lambda);\Z/2)$. Hence the lemma also holds for the degenerate case.
\end{proof}

\section{Proofs of Theorem~\ref{thm_main thm} and corollaries}

Now we are ready to prove the main theorem and its corollaries.

\begin{proof}[Proof of Theorem~\ref{thm_main thm}]
First we prove Statement~\eqref{thm label_Sq^1 on H^2}.
By Proposition~\ref{lemma_cohomology results}~\ref{lemma label_cohmlgy of toric orb}, the number
\[
g=\gcd(d_{i,j}\mid 1\leq i<j\leq n+2)
\]
is the order of $H^3(X(P,\lambda);\Z)$. Since $X(P,\lambda)$ is simply-connected, it has a minimal CW structure~\cite[Theorem~4H.3]{Hat} such that
\[
X(P,\lambda)\simeq\left(\bigvee^{n}_{i=1}S^2\vee P^3(g)\right)\cup e^4.
\]
Since $Sq^1$ acts trivially on $H^*(\bigvee^{n}_{i=1}S^2;\Z/2)$ and the 2- and 3-cells in $P^3(g)$ is connected by mod-$g$ Bockstein homomorphism, $Sq^1\colon H^2(X(P,\lambda);\Z/2)\to H^3(X(P,\lambda);\Z/2)$ is non-trivial if and only if~$g$ is even.

Statement~\eqref{thm label_Sq^1 on H^3} is proved in Lemma~\ref{lemma_Sq^1 on H^3}.

Now we prove Statement~\eqref{thm label_Sq^2}.
Suppose $v_{n+2}$ is a $2$-local smooth vertex.
For $1\leq i\leq n$, let
\[
g_i=\gcd(d_{i,n+1},d_{i,n+2},d_{n+1,n+2}).
\]
By Lemmas~\ref{lemma_Sq on P = Sq on triangle} and~\ref{lemma_criteria of Sq on triangle}, we have
\begin{eqnarray*}
&&\text{$H^*(X(P,\lambda);\Z/2)$ has non-trivial $Sq^2$-action}\\
&\Leftrightarrow&\text{at least one of the $H^*(X(\triangle,\rho_i\lambda);\Z/2)$'s has non-trivial $Sq^2$-action}\\
&\Leftrightarrow&\text{at least one of the $\frac{d_{i,n+1}d_{i,n+2}d_{n+1,n+2}}{g^2_i}$'s is odd}\\
&\Leftrightarrow&
\prod^{n}_{i=1}\left(1-\frac{d_{i,n+1}d_{i,n+2}d_{n+1,n+2}}{g^2_i}\right)\equiv0\pmod{2}.
\end{eqnarray*}
It suffices to show that for $1\leq i\leq n$
\[
\frac{d_{i,n+1}d_{i,n+2}d_{n+1,n+2}}{g^2_i}\equiv\frac{d_{i,n+1} d_{i,n+2}}{g}\pmod{2}.
\]
By assumption $\nu_2(d_{n+1,n+2})=\nu_2(g)=\min\{\nu_2(d_{jk})\mid1\leq j<k\leq n+2\}$. Since
\[
\nu_2(g_i)=\nu_2(\gcd(d_{i,n+1},d_{i,n+2},d_{n+1,n+2}))=\min\{\nu_2(d_{i,n+1}),\nu_2(d_{i,n+2}),\nu_2(d_{n+1,n+2})\},
\]
we have $\nu_2(g_i)=\nu(d_{n+1,n+2})$ and
\[
\frac{d_{i,n+1}d_{i,n+2}d_{n+1,n+2}}{g^2_i}\equiv
\frac{d_{i,n+1}d_{i,n+2}g}{g^2}\equiv
\frac{d_{i,n+1}d_{i,n+2}}{g}\pmod{2}.
\]
Hence the proof is complete.
\end{proof}

\begin{proof}[Proof of Corollary~\ref{coro_spin sufficient}]
By Theorem~\ref{thm_main thm}~\eqref{thm label_Sq^2} it suffices to show that
\[
\prod^{n}_{i=1}\left(1-\frac{d_{i,n+1}d_{i,n+2}}{g}\right)\equiv1\pmod{2},
\]
or equivalently for $1\leq i\leq n$ we have
\[
\frac{d_{i,n+1}d_{i,n+2}}{g}\equiv0\pmod{2}.
\]
Fix an $i\in\{1,\ldots,n\}$ and write $r=\nu_2(g), s=\nu_2(d_{i,n+1})$, and $t=\nu_2(d_{i,n+2})$.
Then
\[
\frac{d_{i,n+1}d_{i,n+2}}{g}\equiv 2^{s+t-r}\pmod{2}.
\]
Since $r$ is the minimum of the three numbers, the index $s+t-r$ is zero if and only if $r=s=t=0$.
However, by assumption $H^*(X(P,\lambda);\Z)$ has $2$-torsion, implying that $g$ is even and $r\geq1$. It follows that $Sq^2$ acts trivially on $H^*(X(P,\lambda);\Z/2)$.
\end{proof}

\begin{proof}[Proof of Corollary~\ref{coro_stable splitting}]
Due to Theorem~\ref{thm_main thm}~\eqref{thm label_Sq^2}, it suffices to show that
\[
\Sigma X(P,\lambda)\simeq\Sigma\C\PP^2\vee\bigvee^{n-1}_{i=1}S^3\vee P^4(g)
\]
if $Sq^2$ acts non-trivially on $H^*(X(P,\lambda);\Z/2)$; otherwise
\[
\Sigma X(P,\lambda)\simeq S^5\vee\bigvee^{n}_{i=1}S^3\vee P^4(g).
\]

Write $g=2^rg'$ where $r=\nu_2(g)$ and $g'$ is odd. By~\cite[Lemma 4.5]{FSS0} we have
\[
X(P,\lambda)\simeq\widehat{X}\vee P^3(g')
\]
where $\widehat{X}=\left(\bigvee^{n}_{i=1}S^2\vee P^3(2^r)\right)\cup e^4$.
Since $H^*(P^3(g');\Z/2)$ has trivial $Sq^2$-action, $Sq^2$ acts non-trivially on~$H^*(X(P,\lambda);\Z/2)$ if and only if it acts non-trivially on $H^*(\widehat{X};\Z/2)$.

If $g$ is odd, then $\widehat{X}$ consists of~$2$-cells and one $4$-cell only. It is a $4$-dimensional Poincar\'{e} duality complex without $2$-torsion in homology.
The lemma follows from~\cite[Theorem 1.1]{ST} or~\cite[Proposition 2.13]{So} in this case.

If $g$ is even, then $Sq^2$ acts trivially on $H^*(X(P,\lambda);\Z/2)$ by Corollary~\ref{coro_spin sufficient}. We need to show
\[
\Sigma\widehat{X}\simeq S^5\vee\bigvee^{n}_{i=1}S^3\vee P^4(2^r).
\]
Let $\widehat{X}_3=\bigvee^{n}_{i=1}S^2\vee P^3(2^r)$ be the $3$-skeleton of $\widehat{X}$ and let $h\colon S^3\to\widehat{X}_3$ be the attaching map of the~$4$-cell.
It suffices to show that $\Sigma h$ is null homotopic.

By Hilton-Milnor Theorem, $h$ can be written as
\[
h\simeq\sum^{n}_{i=1}a_i\eta_i+b\eta_P+ w
\]
where $a_i\in\Z$ and $b\in\Z/2^{r+1}$ are coefficients, $\eta_i\colon S^3\overset{\eta}{\longrightarrow}S^2_i\hookrightarrow X_3$ is the composite of Hopf map and the canonical inclusion, $\eta_P\colon S^3\overset{\eta}{\longrightarrow}S^2\hookrightarrow P^3(g)\hookrightarrow X_3$ is the composite of Hopf map and inclusions,
and $w$ is a sum of Whitehead products. After suspension, $\Sigma w$ vanishes and
\[
\Sigma h\simeq\sum^{n}_{i=1}\bar{a}_i\Sigma\eta_i+\bar{b}\Sigma\eta_P,
\]
where $\bar{a}_i$ and $\bar{b}$ are the reduction images of $a_i$ and $b$ in $\Z/2$.
The maps $\Sigma\eta_i$ and $\Sigma\eta_P$ are detected by~$Sq^2$.
However, $Sq^2$ acts trivially on $H^*(\Sigma\widehat{X};\Z/2)$, so $\bar{a}_i$ and $\bar{b}$ are zero and $\Sigma h$ is null homotopic. Therefore the lemma holds for this case and the proof is complete.
\end{proof}

\begin{proof}[Proof of Corollary~\ref{cor_cohmlgy rigidity}]
For convenience, denote $X(P,\lambda)$ and $X(P',\lambda')$ by $X$ and $X'$ respectively.

When $g$ is odd, this is the main theorem of~\cite{FSS0}.

When $g$ is even, Corollary~\ref{coro_spin sufficient} implies that $Sq^2$ acts trivially on $H^*(X;\Z/2)$. By Corollary~\ref{coro_stable splitting}, it must be the case that $\Sigma X\simeq S^5\vee \bigvee ^n_{i=1}S^3\vee P^4(g)$ where $n$ is the second Betti number of $X$. Similarly, we have~$\Sigma X'\simeq S^5\vee \bigvee ^{n'}_{i=1}S^3\vee P^4(g)$ where $n'$ is the second Betti number of $X'$. By assumption $H^2(X;\Z)\cong H^2(X';\Z)$, so $n=n'$ and $\Sigma X\simeq\Sigma X'$.
\end{proof}

\begin{proof}[Proof of Corollary~\ref{coro_gauge group}]
It suffices to show that if $Sq^2$ acts trivially on $H^*(X(P,\lambda);\Z/2)$ then
\[
\G_k(X(P,\lambda),G)\simeq\G_k(S^4,G)\times\prod^n_{i=1}\Omega^2G\times\Omega^3G\{g\},
\]
where $\Omega^3G\{g\}=Map^*(P^4(g),G)$.
Otherwise
\[
\G_k(X(P,\lambda),G)\simeq\G_k(\C\PP^2,G)\times\prod^{n-1}_{i=1}\Omega^2G\times\Omega^3G\{g\}.
\]
We can prove these two cases using the argument in the proofs of~\cite[Propositions~6.3 and~6.4]{ST}.
\end{proof}

\begin{proof}[Proof of Corollary~\ref{coro_spin criteria}]
The equivalence of (1) and (2) is well known.
For integer $i\geq 0$, the $i$-th Wu class $v_i\in H^i(X(P,\lambda);\Z/2)$ is defined by $Sq^i(x)=x\cup v_i$. It is known that they satisfy Wu Formula:
\[
v_1=w_1,\quad
v_2=w_2+w^2_1,
\quad\cdots
\]
where $w_j$ denotes the $j$-th Stiefel-Whitney class of $X(P,\lambda)$.
Since $X(P,\lambda)$ is orientable, $w_1=v_1=0$. Therefore $v_2=w_2$, implying that $Sq^2$ acts trivially on $H^*(X(P,\lambda);\Z)$ if and only if $w_2=0$.
This fact, together with Theorem~\ref{thm_main thm}~\eqref{thm label_Sq^2}, then implies the corollary.
\end{proof}

\bibliographystyle{amsalpha}

\end{document}